\documentclass[12pt]{article}

\usepackage{graphicx,tikz,color,url}
\usepackage{amsmath,amssymb,latexsym}
\usepackage{cite}
\usetikzlibrary{arrows}

\newtheorem{theorem}{Theorem}[section]

\newtheorem{lemma}[theorem]{Lemma}
\newtheorem{proposition}[theorem]{Proposition}

\newtheorem{corollary}[theorem]{Corollary}

\newcommand{\proof}{\noindent{\bf Proof.\ }}
\newcommand{\qed}{\hfill $\square$ \bigskip}

\newcommand{\smb}{S_{\rm MB}}
\newcommand{\mb}{R_{\rm MB}}

\newcommand{\lc}{{\rm lc}}

\newcommand{\cR}{{\cal R}}
\newcommand{\cS}{{\cal S}}

\newcommand{\cN}{{\cal N}}

\begin{document}
\title{Maker-Breaker resolving game played on lexicographic products of graphs}

\author{ Savitha K S $^{a,}$\thanks{Email: \texttt{savitha@stpauls.ac.in}} 
\and 
Sandi Klav\v{z}ar $^{b,c,d,}$\thanks{Email: \texttt{sandi.klavzar@fmf.uni-lj.si}}
\and 
Tijo James $^{e,}$\thanks{Email: \texttt{tijojames@gmail.com}} 
}
\maketitle

\begin{center}
$^a$ Department of Mathematics, St.\ Paul's  College, Kalamassery, India\\
\medskip
$^b$ Faculty of Mathematics and Physics, University of Ljubljana, Slovenia \\
\medskip

$^{c}$ Institute of Mathematics, Physics and Mechanics, Ljubljana, Slovenia \\

$^d$ Faculty of Natural Sciences and Mathematics, University of Maribor, Slovenia \\
\medskip

$^e$ Department of Mathematics, Pavanatma College, Murickassery, India\\

\end{center}

\begin{abstract} In the Maker-Breaker resolving game, two players named Resolver and Spoiler alternately select unplayed vertices of a given graph $G$. The aim of Resolver is to select all the vertices of some resolving set of $G$, while Spoiler aims to select at least one vertex from every resolving set of $G$. In this paper, this game is investigated on the lexicographic product of graphs. It is proved that if Spoiler has a winning strategy on a graph $H$ no matter who starts the game, or if the first player has a winning strategy on $H$, then Spoiler always has a winning strategy on $G\circ H$. Special attention is paid to lexicographic products in which the second factor is either complete, or a path, or a cycle. For instance, in $G\circ P_{2\ell}$ and in $G\circ C_{2\ell}$, Resolver always wins, while in $G\circ P_{2\ell+1}$ and in $G\circ C_{2\ell+1}$ the same conclusion holds provided $G$ is free from false twins. On the other hand, Spoiler always wins on $G\circ P_5$. In most of the cases, the corresponding Maker-Breaker resolving number is also determined. 
\end{abstract}
\noindent
{\bf Keywords:} Maker-Breaker game; metric dimension; resolving set; Maker-Breaker resolving game; lexicographic product of graphs 

\medskip\noindent
{\bf AMS Subj.\ Class.\ (2020)}: 05C12, 05C57, 05C76

\section{Introduction}
\label{sec:intro}

Let $G=(V(G),E(G))$ be a simple, finite, connected graph with at least two vertices. A set $W\subseteq V(G)$ is a {\em resolving set} of $G$, if for every pair of distinct vertices $x$ and $y$ of $G$, there exists a vertex $z\in W$ such that $d(x,z)\neq d(y,z)$, where $d(u,v)$ denotes the shortest-path distance between the vertices $u$ and $v$. Among all resolving sets of $G$, one with minimum cardinality is a {\em metric basis} and its cardinality is defined to be the {\em metric dimension} of $G$, denoted by $\dim(G)$. As an entry point into this field of research, we recommend the review article~\cite{tillquist-2023} and recent papers~\cite{dolzan-2025, knor-2025, tapendra-2025}. 

The Maker-Breaker resolving game was introduced in~\cite{kang-2021} and the authors investigated the same in the corona product of graphs~\cite{tijo-2025}. A closely related Maker-{B}reaker strong resolving game has been investigated in~\cite{kang-2024+}, the fractional version of the game in~\cite{yi-2021}, while for the general framework on Maker-Breaker games see Chapter 2 of the book~\cite{hefetz-2014}. For more on these games, see also the recent articles~\cite{bennett-2025, duchene-2025} and the references therein.  

In the  {\em Maker-Breaker resolving game} ({\em MBRG} for short), two players, named Resolver and Spoiler, alternately select unplayed vertices of a given graph $G$. The aim of Resolver is to select all the vertices of some resolving set of $G$, while Spoiler aims to select at least one vertex from every resolving  set of $G$. The player who achieves their goal is the winner of the game. If Resolver starts the game, we speak of an  {\em R-game}, otherwise we speak of an {\em S-game}. 

For the {\em outcome} $o(G)$ of the MBRG  played on $G$ we have $o(G)\in \{\mathcal{R}, \mathcal{S}, \mathcal{N}\}$, with the following meaning \cite{kang-2021}: 
     \begin{itemize}
         \item $o(G)=\mathcal{R}$: Resolver has a winning strategy no matter who starts the game; 
         \item $o(G)=\mathcal{S}$: Spoiler has a winning strategy no matter who starts the game;
         \item  $o(G)=\mathcal{N}$: the first player has a winning strategy.
    \end{itemize}

In addition to knowing who wins the game, it is also of interest how fast the winner can achieve this. The {\em Maker-Breaker resolving number}, $\mb(G)$, is the minimum number of moves of Resolver to win the R-game in $G$ when both players play optimally. The {\em Maker-Breaker spoiling number}, $\smb(G)$, is the minimum number of moves of Spoiler to win the R game in $G$ when both players play optimally. For the S-game, the corresponding invariants are denoted by $\mb'(G)$ and $\smb'(G)$. 

The lexicographic product $G\circ H$ of graphs $G$ and $H$ has the vertex set $V(G\circ H)=\{(g,h):\ g\in V(G), h\in V(H)\}$ and the edge set $E(G\circ H)= \{(g,h)(g',h'):\ gg' \in E(G), ~\text{or}~ g=g' ~\text{and}~ hh'\in E(H)\}$.  Studies on lexicographic products in various domains highlight its versatility, which offers valuable insights into network structures, information flow, and resilience. In \cite{jann-2012}, the metric dimension of the lexicographic products is determined in terms of a parameter called the adjacency dimension. The general bounds for the metric dimension of the lexicographic product of graphs can be seen in \cite{saputro-2013}. For some other studies of the lexicographic product of graphs we refer to ~\cite{bal-2025, cabrera-2024, jann-2025, kuziak-2023, orel-2023, samod-2021}. 

In this paper, we investigate the MBRG played on the lexicographic product of graphs. The following section is preparatory and, among other things, recalls known results that we need here. In Section 3, we discuss some general properties of the MBRG on lexicographic products. In particular, if $G$ and $H$ are nontrivial connected graphs and $o(H)\in\{\cN,\cS\}$, then we show that $o(G\circ H)=\mathcal{S}$. In Section~\ref{sec:complete} we consider lexicographic products in which one factor is complete. We prove that if $G$ has true twins, then $\smb(G\circ K_2) = \smb'(G\circ K_2) = 2$, otherwise $\mb(G\circ K_2) = \mb'(G\circ K_2) = n(G)$. In Section~\ref{sec:paths} we focus on the MBRG where the second factor is a path or a cycle, while in the concluding section we pose some open problems. 

\section{Preliminaries}
\label{sec:prelim}

In this section, we collect known results and concepts which will be used in the rest of the paper, and add some additional definitions. Let's start with the latter. 

For a positive integer $\ell$, we will denote the set $\{1,\ldots, \ell\}$ by $[\ell]$. By $n(G)$ we denote the order of $G$. For a vertex $v\in V(G)$, its open neighborhood is $N_G(v)=\{u\in V(G):\ uv\in E(G)\}$. Vertices $u$ and $v$ are {\em twins} if $N(u)\setminus\{v\} = N(v)\setminus \{u\}$. The twin relation is an equivalence relation, and an equivalence class under this equivalence relation is called a {\em twin equivalence class}. A pair of adjacent vertices in a twin equivalence class is termed as {\em true twins}, and a pair of non-adjacent vertices in a twin equivalence class is termed as {\em false twins}.

A set $S\subseteq V(G)$ is a {\em dominating set} of a graph $G$, if every vertex $u\in V(G)\setminus S$ has a neighbor in $S$. A vertex $x\in V(G)$ is a {\em dominating vertex} of $G$ if $x$ is adjacent to all the other vertices of $G$. Note that in this case $\{x\}$ is a dominating set of $G$.

If $g\in V(G)$, then by $^gH$ we denote the {\em $H$-layer} of $G\circ H$ corresponding to $g$, that is, the subgraph of $G\circ H$ induced by the vertices of the form $(g,h)$, where $h$ runs over vertices of $H$. For $h\in V(H)$, the {\em $G$-layer} $G^h$ is defined analogously. Each $G$-layer $G^h$ is isomorphic to $G$, and each $H$-layer $^gH$ is isomorphic to $H$.

Let $A=\{\{u_1, v_1\},\ldots, \{u_k, v_k\}\}$ be a set of 2-subsets of $V(G)$ of pairwise different elements, that is, $|\bigcup_{i=1}^k{\{u_i,v_i\}}|=2k$. Then $A$ is called a {\em pairing resolving set} if every set $\{x_1,\ldots, x_k\}$, where $x_i\in \{u_i,v_i\}$, is a resolving set of $G$. A pairing resolving set $A$ of $G$ with $|A| = \dim(G)$ is a {\em dim-pairing resolving set} of $G$. 

$S\subseteq V(G)$ is a {\em locating set} of $G$ if $N_G(u)\cap S \ne N_G(v)\cap S$ for every two vertices $u,v\in V(G)\setminus S$. By $\lc(G)$ we denote the size of a smallest locating set of $G$. If in addition $N_G(u)\cap S \ne S$ for every $u\in V(G)\setminus S$, then $S$ is a {\em strictly locating set} of $G$. 

In the second part of this section, we recall the known results that we need is the following. 

\begin{proposition} {\rm \cite[Proposition 2.3]{kang-2021}}
\label{prop:bound dimension} 
If $G$ is a connected graph, then the following properties hold.
\begin{enumerate}
    \item[(i)] If $o(G)=\cR$, then $\mb'(G)\geq\mb(G)\geq \dim(G)$.
    \item[(ii)] If $o(G)=\cS$, then $\smb(G)\geq\smb'(G)$.
\end{enumerate}
\end{proposition}

\begin{proposition} {\rm \cite[Proposition~3.3]{kang-2021}}
\label{prop:pairing} 
  If $G$ admits a pairing resolving set, then $o(G) =\mathcal{R}$.
  Moreover,  if $G$ admits a dim-pairing resolving set, then $\mb(G) = \mb'(G)
 = \dim(G)$.
\end{proposition}
  
\begin{proposition}{\rm \cite[Proposition~3.2]{kang-2021}}
 \label{prop:twin equivalence} 
 Let $G$ be a connected graph with $n(G)\geq 4$.
 \begin{itemize}
     \item[(a)]  If $G$ has a twin equivalence class of cardinality at least $4$, then $o(G)=\mathcal{S}$ and $\smb(G) = \smb'(G)=2$.
     \item[(b)]  If $G$ has two twin equivalence classes, each of cardinality at least $3$, then $o(G)=\mathcal{S}$ and $\smb(G) = \smb'(G)=2$.
 \end{itemize} 
\end{proposition}

\begin{theorem} {\rm\cite{monsanto-2022}}
\label{thm:locating}
Let $G$ and $H$ be nontrivial connected graphs, where $\Delta(H) \leq n(H)-2$. Then $W =\cup_{x\in S}[\{x\} \times T_x]$, where $S\subseteq V(G) $ and $T_x \subseteq V(H)$ for
each $x \in S$, is a resolving set of $G\circ H$ if and only if
\begin{itemize}
    \item(i) $S = V(G)$;
\item(ii) $T_x$ is a locating set of $H$ for every $x \in V(G)$;
\item (iii) $T_x$ or $T_y$ is a strictly locating set of $H$ whenever $x$ and $y$ are adjacent
vertices of $G$ with $N_G[x] = N_G[y]$.
\item (iv) $T_x$ or $T_y$ is a (locating) dominating set of $H$ whenever $x$ and $y$ are nonadjacent vertices of $G$ with $N_G(x) = N_G(y)$.
\end{itemize}
\end{theorem}

\begin{lemma}{\rm\cite{tijo-2025}}
\label{lem:path-strictly-locating}   Let $\ell \ge 3$ and let $V(P_{2\ell}) = V(C_{2\ell}) = [2\ell]$. If $W=\{w_{1},\ldots, w_{\ell}\}$, where $w_{i}\in \{2i-1, 2i\}$ for $i\in [\ell]$, then $W$ is a strictly locating set of $P_{2\ell}$ as well as a strictly locating set of $C_{2\ell}$.
\end{lemma}

\section{Some general properties}
\label{sec:general}

In this section, we present some general properties of the MBRG played on $G\circ H$. We first observe that if $o(G\circ H)=\cR$, then 
\begin{equation}
n(G)\dim(H)\leq \mb(G\circ H) \leq \mb'(G\circ H) \leq \left\lfloor\frac{n(G)n(H)}{2}\right \rfloor. \label{eq:1} 
\end{equation}

The upper bound is obvious because the players alternate their moves. The lower bound follows from the inequality  $n(G)\dim(H)\leq \dim(G\circ H)$ established in~\cite{saputro-2013}. 

By Theorem~\ref{thm:locating}, every resolving set of $G\circ H$ intersects each of its $H$-layers in a locating set of $H$. Together with \eqref{eq:1} this yields: 

\begin{corollary}
\label{cor:bound}
 If $o(G\circ H)=\cR$ and $\Delta(H)\leq n(H)-2$, then $\mb(G\circ H) \ge n(G)\lc(H)$.
\end{corollary}

Consider the lexicographic product $P_4\circ P_4$. The strategy of Resolver to select two vertices in each $^uP_4$-layer is a winning one no matter who starts the game, hence $o(P_4\circ P_4) = \cR$. Since $\lc(P_4) = 2$, Corollary~\ref{cor:bound} implies $\mb'(P_4\circ P_4) \ge \mb(P_4\circ P_4) \ge 8$. Then by the upper bound of~\eqref{eq:1} we get $\mb(P_4\circ P_4) = \mb'(P_4\circ P_4) = 8$. This demonstrates that the upper bound of~\eqref{eq:1} is sharp. Moreover, since $\mb(P_4\circ P_2) = \mb'(P_4\circ P_2) = 4$, the lower bound is sharp as well. In fact, in the latter example, the two bounds coincide.

Concerning the outcome, we have the following general result for all the cases when $o(H)\in\{\cN,\cS\}$. 

\begin{theorem}
\label{thm:outcome}
If $G$ and $H$ are nontrivial connected graphs and $o(H)\in\{\cN,\cS\}$, then $o(G\circ H)=\mathcal{S}$. 
\end{theorem}

\proof
Let $S$ be a resolving set of $G\circ H$. For $g\in V(G)$, let $S_g=S\cap V(^gH)$. Since every two vertices of $V(^gH)$ have the same distance to any vertex from $V(G\circ H)\setminus V(^gH)$, the set $S_g$ is a resolving set of $H$. 

The strategy of Spoiler is to play the first move in some $^gH$-layer in which Resolver has not yet played. As $G$ is nontrivial, Resolver can achieve this goal irrespective of who starts the game. After that, since $o(H)=\{\cN,\cS\}$, Spoiler's strategy is to play only in $^gH$  by following her optimal strategy in $H$. In this way, she wins the game in $^gH$ and thus also in $G\circ H$, because in view of the No-Skip-Lemma~\cite[Lemma 2.2]{kang-2021}, her goal is not spoiled if she is forced to play some consecutive moves in $^gH$. Thus $o(G\circ H)=\mathcal{S}$.
\qed

As we will see in the rest of the paper, no theorem parallel to Theorem~\ref{thm:outcome} can be stated when $o(H) = \mathcal{R}$, cf.\ for instance, Theorem~\ref{thm:K2}.  

\section{The game on products with a complete factor}
\label{sec:complete}

In this section, we consider lexicographic products where one of the factors is complete. The case when both factors are complete is trivial, since in this case, the graph in question is complete.  

Let $G$ be a nontrivial connected graph. Since $o(K_3) = \cN$ and $o(K_n) = \cS$ for every $n\ge 4$, Theorem~\ref{thm:outcome} implies that $o(G\circ K_n) = \mathcal{S}$, $n\geq 3$. Moreover, using the general strategy from the proof of Theorem~\ref{thm:outcome}, Spoiler can finish the game in two moves. Hence: 

\begin{proposition}
If $G$ is a nontrivial connected graph and $n\ge 3$, then $\smb(G\circ K_n) = \smb'(G\circ K_n) = 2$. 
\end{proposition}

With respect to complete second factors, the remaining case is $K_2$, for which we have $o(K_2) = \mathcal{R}$. 

\begin{theorem}
\label{thm:K2}
If $G$ is a nontrivial connected graph, then $o(G\circ K_2)\in \{\cR, \cS\}$. Moreover, the following holds.  
\begin{itemize}
    \item If $G$ has true twins, then $\smb(G\circ K_2) = \smb'(G\circ K_2) = 2$. 
    \item If $G$ is free from true twins, then $\mb(G\circ K_2) = \mb'(G\circ K_2) = n(G)$.
\end{itemize}
\end{theorem}

\proof 
Let $G$ be a nontrivial connected graph with $V(G)=\{g_1, \ldots, g_{n}$\}, and let $V(K_2) = [2]$. 

Assume first that $G$ contains true twins and let $g_i$, $g_j$, $i\neq j$, be such vertices. Then the vertices $(g_i, 1)$, $(g_i, 2)$, $(g_j, 1)$, $(g_j, 2)\}$ belong to the same twin equivalence class of $G\circ K_2$. Therefore, $G\circ K_2$ has a twin equivalence class of cardinality at least $4$. By Proposition~\ref{prop:twin equivalence} we get $ o(G\circ K_2)=\mathcal{S}$, as well as $\smb(G\circ K_2) = \smb'(G\circ K_2) = 2$. 

Assume second that $G$ has no true twins. Then we claim that the set $S = \{(g_i, k_i):\ i\in [n], k_i\in \{1,2\}\}$ forms a resolving set of $G\circ K_2$. Consider arbitrary vertices which are not in $S$, say $(g_i, k'_{i})$ and $(g_j,k'_j)$, where $i\neq j$, $k'_i = 3 - k_i$, and $k'_j = 3 - k_j$. If $g_ig_j\notin E(G)$, then $d((g_i, k'_{i}), (g_i, k_i)) = 1$ and $d((g_j, k'_{j}), (g_i, k_i)) > 1$. Similarly, if $g_ig_j\in E(G)$, then since $g_i$ and $g_j$ are not true twins, there exists a vertex $g_l$ such that $g_ig_l\in E(G)$ and $g_jg_l\notin E(G)$. Now $d((g_i, k'_{i}), (g_l, k_{l})) = 1$ and $d((g_j, k'_{j}), (g_l, k_{l}) > 1$. 

The above argument implies that $\{\{(g_i, 1), (g_i, 2)\}: i\in [n]\}$ is  a pairing resolving set of $G \circ K_2$. Since by Theorem~\ref{thm:locating} this set is also a dim-pairing resolving set, the required conclusion when $G$ has no true twins follows by Proposition~\ref{prop:pairing}.
\qed

When the first factor of a lexicographic product is complete, we have the following partial result. 

\begin{proposition}
If $G$ is a nontrivial connected graph with $n(G)\geq 4$, then the following hold.
\begin{enumerate}
\item[(i)] If $G$ has a dominating vertex, then $o(K_m\circ G)=\cS$ for $m\geq 4$.
\item[(ii)] If $G$ has two dominating vertices, then $o(K_m\circ G)=\cS$ for $m\in \{2,3\}$.
\end{enumerate}
Moreover, in each of the cases, Spoiler can win the game in two moves. 
\end{proposition}

\proof
(i) Let $g$ be a dominating vertex of $G$. Then each vertex from $K_m^g$ is a dominating vertex in $K_m\circ G$. Hence, these vertices form a twin equivalence class. Since $n(G)\geq 4$, the assertion follows by Proposition~\ref{prop:twin equivalence}. 

(ii) Let $g$ and $g'$ be two dominating vertices of $G$. In this, each vertex from the set $V(K_m^g)\cup V(K_m^{g'})$ is a dominating vertex of $K_m\circ G$, hence we can again apply Proposition~\ref{prop:twin equivalence} to reach the desired conclusion. 
\qed

\section{MBRG where the second factor is a path or a cycle}
\label{sec:paths}

In this section, we consider the MBRG on $G\circ P_n$ and on $G\circ C_n$. We assume throughout the section that $V(G) = \{(g_i:\ i\in [n(G)]\}$ and $V(P_n) = V(C_n) = [n]$, with the natural adjacencies, and in cycles with indices computed modulo $n$ when needed. 

\subsection{Short paths}

We first consider the case of paths $P_n$, where $n\in \{3,4,5\}$.  

\begin{theorem}
If $G$ is a connected twin-free graph then $o(G\circ P_3)=\cR$ and $\mb(G\circ P_3) = \mb'(G\circ P_3) = n(G)$.
\end{theorem}

\proof
Let $R=\{\{(g_i,1),(g_i,3)\}:\ i\in [n(G)]\}$. We claim that $R$ is a dim-pairing resolving set of $G\circ P_3$. 

Let $s_i$, $i\in [n(G)]$ be an arbitrary, fixed integer from the set $\{1,3\}$, and let $S=\{(g_i,s_i): i\in [n(G)]\}$. We first need to show that $S$ is a resolving set of $G\circ P_3$. For this, first consider the vertices of a $^{g_i} P_3$-layer. Since $d((g_i,1),(g_i,2)) = 1 = d((g_i,2),(g_i,3))$ and $d((g_i, 1),(g_i, 3))=2$, each vertex in this layer is resolved by $(g_i, s_i)$. Next, consider two  vertices from different $P_3$-layers, say $(g_\ell, j)$ and $(g_k, j')$, where $j,j' \in [3]$. Since $G$ is twin-free, there is a vertex $g_h$ in $V(G)$ that is adjacent to $g_\ell$ but not to $g_k$. Then $d((g_\ell, j),(g_h, s_h))= 1$ and $d((g_k,j'),(g_h, s_i)) \neq 1$, therefore, these two vertices are resolved by $(g_h, s_i) \in S$. We have thus proved that $R$ forms a pairing resolving set. Moreover, by Theorem~\ref{thm:locating} we get that $\dim(G\circ P_3) = n(G)$, hence $R$ is indeed a dim-pairing resolving set as claimed. In view of Proposition~\ref{prop:pairing} we can conclude that  $\mb(G\circ P_3) = \mb'(G\circ P_3) = n(G)$. 
\qed

\begin{theorem}
If $G$ is a nontrivial connected graph, then $o(G\circ P_4)=\cR$ and $\mb(G\circ P_4) = \mb'(G\circ P_4) = 2 n(G)$. 
\end{theorem}

\proof
As per the characterization of resolving sets of lexicographic products in Theorem~\ref{thm:locating}, we describe a strategy for Resolver which allows him to select the appropriate set as follows. 

Let $g_i$ and $g_j$ be false twins of $G$. Assume that Spoiler selects $(g_i, 1)$ as her first move in $^{g_i}P_4$. Then the aim of Resolver is to form a locating dominating set in $^{g_i}P_4$. Thus, Resolver selects $(g_i, 2)$ as his optimal response. Clearly $\{(g_i, 2), (g_i, 3)\}$ and $\{(g_i, 2), (g_i, 4)\}$ are locating dominating sets. Therefore, Resolver can form a locating dominating set in $^{g_i}P_4$ after any move of Spoiler in this layer. If Spoiler selects $(g_i, 2)$ as her first move in $^{g_i}P_4$, then Resolver selects $(g_i, 1)$ as his optimal response. Then, Resolver will form a locating dominating set in $^{g_i}P_4$ after the next move of Spoiler in $^{g_i}P_4$, by selecting either $(g_i, 3)$ or $(g_i, 4)$. Other moves of Spoiler in $^{g_i}P_4$ are dealt in a similar fashion by Resolver.

Next, consider the case when $g_i$ and $g_j$ are true twins of $G$. Assume that Spoiler selects $(g_i, 1)$ as her first move in $^{g_i}P_4$. Here, the aim of Resolver is to form a strictly locating set. Therefore, Resolver selects $(g_i, 3)$  as his next optimal response. Now $\{(g_i, 2), (g_i, 3)\}$  and $\{(g_i, 3), (g_i, 4)\}$  are strictly locating sets of $^{g_i}P_4$. Thus, by his second move, Resolver succeeds in forming a strictly locating set in $^{g_i}P_4$. If Spoiler selects $(g_i, 2)$ as her first move, then Resolver selects $(g_i, 4)$ as his optimal response. Clearly, Resolver can form a strictly locating set by choosing either $(g_i, 1)$ or $(g_i, 3)$, depending on Spoiler's move. Similar arguments hold for other cases.

Let $g_i$ and $g_j$ be any two vertices of $G$ other than the above-mentioned cases. Then the fact that any pair of vertices in $P_4$ forms a locating set helps Resolver to obtain a locating set in each $P_4$ layer, whatever the Spoiler's move in that layer. 
Hence, in all cases, Resolver wins the game on $G\circ P_4$ as a second player. Thus $o(G\circ P_4)=\cR.$

The equalities $\mb(G\circ P_4)=\mb'(G\circ P_4) = 2n(G)$ follow from Corollary~\ref{cor:bound} because $\lc(P_4)=2$.
 \qed

\begin{theorem}
If $G$ is a nontrivial connected graph, then
$o(G\circ P_5)=\cS$. Moreover $\smb'(G\circ P_5)=\smb(G\circ P_5)=3$. 
\end{theorem}

\proof
Consider the $R$-game on $G\circ P_5$ and the following strategy of Spoiler. 

If Resolver selects a vertex in $^{g_i}P_5$, then Spoiler selects the vertex $(g_j, 3)$, where $i\neq j$. After that, Spoiler's strategy is the following. If  Resolver selects a vertex from $\{(g_j, 1), (g_j, 5)\}$, then Spoiler replies by selecting the other vertex from the set. If Resolver selects a vertex from $\{(g_j, 2), (g_j, 4)\}$, then again Spoiler selects the other vertex from the set. In this way, in $^{g_j}P_5$, Resolver will have one of the following choices: $\{(g_j, 1), (g_j, 2)\}$, $\{(g_j, 1), (g_j, 4)\}$, $\{(g_j, 2), (g_j, 5)\}$, $\{(g_j, 4), (g_j, 5)\}$.
Since each set in the above collection fails to form a locating set in $^{g_j}P_5$, Theorem~\ref{thm:locating} implies that Resolver fails to select a resolving set in $G\circ P_5$. Hence, Spoiler wins the game on $G\circ P_5$ as the second player, and therefore she also wins the game as the first player. So $o(G\circ P_5)=\cS$. Moreover $\smb'(G\circ P_5)\leq \smb(G\circ P_5)\leq 3$. Also, it is clear that any three vertices from a $P_5$-layer always form a strictly locating set of the layer.  Therefore, Spoiler cannot win in $G\circ P_5$ in two moves. Thus $ \smb'(G\circ P_5)\geq 3$ and we can conclude that $\smb'(G\circ P_5) = \smb(G\circ P_5) = 3$.
\qed 

\subsection{Short cycles}
We next consider the case of cycles $C_n$, where $n\in \{4,5\}$. (The case $n=3$ is covered in Section~\ref{sec:complete}.)  

\begin{theorem}~\label{thm:c_4}
If $G$ is a nontrivial connected graph, and $n\in \{4, 5\}$, then
$o(G\circ C_n) =\cR$ and $\mb(G\circ C_n)=\mb'(G\circ C_n) = 2 n(G)$.
\end{theorem}

\proof 
Each pair of adjacent vertices of $C_4$ forms a strictly locating and dominating set of $C_4$. Hence, in view of Theorem~\ref{thm:locating}, if Resolver can select two adjacent vertices in each $C_4$-layer, then he wins. This is indeed not difficult for him to achieve. If the first vertex selected in $^{g_i}C_4$ is by Spoiler, and this vertex is $(g_i,j)$, then Resolver responds by selecting $(g_i, j+2)$. Afterwards, Resolver is able to select one vertex from $\{(g_i, j+1), (g_i, j+3)\}$, thereby achieving his goal. Resolver can follow this strategy in each $C_4$-layer, hence he wins as the second player in $G\circ C_4$ and thus also as the first player. Therefore $o(G\circ C_4)=\mathcal{R}$. 

The case $G\circ C_5$ can be treated similarly. More precisely, note first that any pair of adjacent vertices of $C_5$ forms a strictly locating set of $C_5$, and any pair of nonadjacent vertices in $C_5$ forms a locating dominating set of $C_5$. In view of Theorem~\ref{thm:locating}, a strategy of Resolver is to form a strictly locating set of $C_5$ in $^{g_i}C_5$ when $g_i$ is a true twin vertex, a locating dominating set of $C_5$ in $^{g_i}C_5$ when $g_i$ is a false twin vertex, and a resolving set in any other $C_5$-layer. If Spoiler selected $(g_i,j)$ as the first vertex played in $^{g_i}C_5$, and $g_i$ is a true twin vertex, then Resolver responds by selecting $(g_i, j+2)$. In the other two cases, Resolver responds by selecting $(g_i, j+1)$. After that, in each of the cases, Resolver can achieve his goal. Hence $o(G\circ C_5)=\mathcal{R}$. 

Since  $\dim(C_4) = 2$, by~\eqref{eq:1} we get 
$$2n(G)\leq \mb(G\circ C_4)\leq \mb'(G\circ C_4)\leq \left\lfloor\frac{n(G) n(C_4)}{2}\right \rfloor=2 n(G)\,.$$ 
Similarly, since $\dim(C_5) = 2$, we get by~\eqref{eq:1} that $\mb(G\circ C_5) \ge 2n(G)$, and we can conclude, having in mind the above strategy of Resolver, that also in this case  $\mb(G\circ C_5)=\mb'(G\circ C_5) = 2 n(G)$. 
\qed

\subsection{Longer paths and cycles}

We now focus on paths and cycles of length at least $6$. 

\begin{theorem}\label{thm:path on 2l vertices}
If $G$  is a connected nontrivial  graph, then $o(G\circ P_{2\ell})  = o(G\circ C_{2\ell}) = \cR$ for $\ell\geq 3$. Moreover,  $\mb(G\circ P_6)= \mb'(G\circ P_6)=3n(G)$. 
\end{theorem}

\proof
We first consider the paths $P_{2\ell}$, $\ell\geq 3$. In order to show that $o(G\circ P_{2\ell}) = \cR$, Resolver follows the strategy that in  $^{g_i}P_{2\ell}$ he selects one vertex from each of the sets $\{(g_i,2j-1), (g_i,2j) \}$, $j\in [\ell]$. This strategy is clearly feasible regardless of who starts the game. Using it, Resolver selects in $^{g_i} P_{2\ell}$ a set $W_{g_i}$ which is by Lemma~\ref{lem:path-strictly-locating} a strictly locating set of $P_{2\ell}$. Since $W_{g_i}$ consists of one of the end vertices of each edge of the perfect matching of $P_{2\ell}$, the set $W_{g_i}$ is also a dominating set of $P_{2\ell}$. As this holds for any $P_{2\ell}$-layer, by Theorem~\ref{thm:locating} we can conclude that $o(G\circ P_{2\ell})=\cR$.

The argument that $o(G\circ C_{2\ell})=\cR$ holds for $\ell\geq 3$ is completely parallel to the above argument for paths, hence we do not repeat it here. 

It remains to demonstrate that $\mb(G\circ P_6)= \mb'(G\circ P_6) = 3n(G)$. For this sake note first that $\lc(P_6) = 2$ and that $\{2, 4\}$ and $\{3, 5\}$ are the only locating sets of $P_6$ with exactly two vertices. In the MBRG game played on $G\circ P_6$, if Resolver selects  $(g_i,2)$, then Spoiler selects $(g_i,4)$ and vice versa. In a similar manner, if  Resolver selects a vertex from $\{(g_i,v_3), (g_i, v_5)\}$, then Spoiler selects the other one from the set. Thus, in each $P_6$-layer of $G\circ P_6$, Spoiler has a strategy which ensures that Resolver needs at least three moves to form a locating set of $P_6$. This strategy is independent of who first selects a vertex in a $P_6$-layer. Thus $\mb'(G\circ P_6)\geq \mb(G\circ P_6)\geq 3n(G)$. On the other hand, by Theorem~\ref{thm:path on 2l vertices} we get $\mb(G\circ P_6)\leq \mb'(G\circ P_6)\leq 3n(G)$. 
\qed

\begin{theorem}
\label{thm:long}
If $G$  is a connected nontrivial graph that is free from false twins, then $o(G\circ P_{2\ell+1}) = o(G\circ C_{2\ell + 1}) = \cR$ for $\ell\geq 3$.
\end{theorem}

 \proof
 First, we prove this for paths $P_{2\ell+1}$, $\ell\geq 3$. For any $g_i\in V(G)$, partition the vertex set of $^{g_i}P_{2\ell+1}$ as follows: 
 $$\Big\{  \{(g_i, 1), (g_i, 2)\}, \ldots, 
             \{(g_i, 2\ell-3), (g_i, 2\ell-2)\},\{(g_i, 2\ell-1), (g_i, 2\ell), (g_i, 2\ell+1)\} \Big\}\,.$$
The strategy of Resolver is the following. As soon as Spoiler selects a vertex from  $\{(g_i, 2j-1), (g_i, 2j\}$, $j\in[\ell-1]$, Resolver selects the other vertex from the set. If Spoiler selects a vertex from  $\{ (g_i, 2\ell-1), (g_i, 2\ell), (g_i, 2\ell+1)\}$, then we distinguish the following cases. 

\medskip\noindent
{\bf Case 1}: Spoiler has selected $(g_i,2\ell-1)$. \\
In the subcase when $(g_i, 2\ell-2)$ has not yet been selected, Resolver plays $(g_i, 2\ell-2)$. Then, before the end of the game, Resolver will be able to select one of $(g_i,2\ell)$ and $(g_i,2\ell+1)$, hence in this way a strictly locating set will be selected by him in $^{g_i}P_{2\ell+1}$. Consider next the subcase when $(g_i, 2\ell-2)$ has already been selected. If $(g_i,2\ell-2)$ has been selected by Resolver (and hence $(g_i, 2\ell-3)$ by Spoiler), then Resolver selects next the vertex $(g_i, 2\ell)$. And if $(g_i, 2\ell-2)$ has been selected by Spoiler (and hence $(g_i, 2\ell-3)$ by Resolver), then Resolver next selects $(g_i, 2\ell+1)$. In each of the cases, the vertices selected by Resolver form a strictly locating set of $^{g_i}P_{2\ell+1}$ layer.

\medskip\noindent
{\bf Case 1}: Spoiler has selected $(g_i, 2\ell)$ or $(g_i, 2\ell+1)$. \\
In this case, Resolver replies by selecting $(g_i, 2\ell-1)$. We can then argue as above that in this way, Resolver has constructed a strictly locating set of $^{g_i}P_{2\ell+1}$. 

\medskip
The above strategy of Resolver is performed in each of the $P_{2\ell+1}$-layers, hence by Theorem~\ref{thm:locating} we can conclude $o(G\circ P_{2\ell+1})=\cR$.

\smallskip
Consider next cycles $C_{2\ell+1}$, $\ell\geq 3$.  We may, without loss of generality, assume that Spoiler starts the game by selecting the vertex $(g_i, 2\ell+1)$. Resolver replies by playing the vertex $(g_i,2\ell)$. After that, Resolver continues using the following strategy on $^{g_i}C_{2\ell +1}$. Consider the subsets 
$$Z_{g_i, j} = \{(g_i, 2j-1), (g_i,2j) \}, j\in [\ell-2]\,,$$ 
and set also
$$Z_{g_i, \ell-1} = \{(g_i,2\ell-3), (g_i,2\ell-2), (g_i,2\ell-1)\}\,.$$ Whenever Spoiler selects a vertex from some set $Z_{g_i, j}$, $j\in [\ell-2]$, Resolver replies by playing the other vertex from the set. At some point, Spoiler will select a vertex from $Z_{g_i, \ell-1}$. If this vertex is $(g_i,2\ell-3)$ or $(g_i, 2\ell-1)$, then Resolver replies by playing $(g_i, 2\ell-2)$.  And if the first vertex from $Z_{g_i, \ell-1}$ selected by Spoiler is $(g_i, 2\ell-2)$, then Resolver replies by playing $(g_i, 2\ell-1)$.

We claim that using the above-described strategy, Resolver constructs a strictly locating set of $^{g_i}C_{2\ell +1}$. To show it, we need to demonstrate that no four consecutive vertices were selected by Spoiler and that there are no five consecutive vertices such that only the middle one was selected by Resolver. The first situation cannot happen because in every set $Z_{g_i,j}$ Resolver selected one vertex and because in $Z_{g_i, \ell-1}$ he has also selected one vertex. The second situation could only happen if Spoiler selects $(g_i, 2\ell+1)$ and $(g_i,1)$, but then since Resolver selected one of the vertices $(g_i, 2\ell-1)$ and $(g_i, 2\ell-2)$, this also does not happen in this case. This proves the claim. 

By the above argument and using the No-Skip Lemma~\cite[Lemma 2.2]{kang-2021}, Resolver can select a strictly locating set in each $C_{2\ell +1}$-layer. Therefore, by Theorem~\ref{thm:locating} we can conclude that $o(G\circ C_{2\ell + 1}) = \cR$.  
\qed

The assumption of Theorem~\ref{thm:long} that $G$ is free from false twins cannot be avoided. For this sake, consider $G=K_{1,3}$ with $V(G)=\{g_0,g_1,g_2,g_3\}$, where $g_0$ is the vertex of degree $3$. Clearly, the vertices $g_1, g_2$ and $g_3$ are false twins. Consider now the lexicographic product $G\circ P_7$ and let Spoiler first select the vertex $(g_1, 2)$. The aim of Resolver is to form a locating dominating set in $^{g_1}P_7$. For this reason, Resolver must select the vertex $(g_1, 1)$. Afterwards, Spoiler selects $(g_1, 6)$, which in turn forces Resolver to select $(g_1,7)$. As the third move, Spoiler selects $(g_1, 4)$, which enables him to later select one of $(g_1,3)$ and $(g_1,5)$, depending on the next move of Resolver. In this way, Resolver fails to select a locating-dominating set in $^{g_1}P_7$. Since Spoiler can further apply the above-described strategy also in at least one of $^{g_2}P_7$ and $^{g_3}P_7$, by the No-Skip Lemma and Theorem~\ref{thm:locating} we conclude that $o(G\circ P_7)\neq \cR$.

\section{Concluding remarks}
\label{sec:conclude}

We conclude the paper with the following open problems.
\begin{itemize} 
\item Determine $o(G\circ P_{2\ell+1})$ and $o(G\circ C_{2\ell + 1})$ when $G$ contains false twins.

\item Characterize the graphs $G$ and $H$ for which $o(H)=\cR$ and $o(G \circ H)=\cR$.

\item If $G$ is a nontrivial connected graph without true twins, then $\mb(G\circ P_2)=n(G)$, see Theorem~\ref{thm:K2}. Further, by Theorem~\ref{thm:c_4} and Theorem~\ref{thm:path on 2l vertices}, $\mb(G\circ P_4)=2n(G)$ and $\mb(G\circ P_6)=3n(G)$. It remains to determine $\mb(G\circ P_{2\ell})$ for $\ell\geq 4$.
\end{itemize}

\section*{Declaration of interests}
 
The authors declare that they have no conflict of interest. 

\section*{Data availability}
 
Our manuscript has no associated data.

\section*{Acknowledgements}

Sandi Klav\v{z}ar was supported by the Slovenian Research and Innovation Agency (ARIS) under the grants  P1-0297, N1-0285, N1-0355, N1-0431.

\end{document}